\documentclass{article}
\title{Stasis Points and Approximating Two-Cycles}
\author{Stewart D. Johnson\\Department of Mathematics and Statistics\\Williams College}

\newcommand{\qed}{\vbox{\hrule width4pt height3pt depth1pt} }

\newcommand{\emp}[1]{\bf #1 \rm}

\font\bbb=msbm10
\newcommand{\R}{\hbox{\bbb\char82}}

\newcommand{\del}{\partial}

\newcommand{\ident}{\equiv}

\newtheorem{theorem}{Theorem}

\newenvironment{proof}{{\bf Proof: }}{\hfill$\qed$}

\begin{document}             % End of preamble and beginning of text.
%\openup 8pt
\maketitle

\begin{abstract}
\noindent Generically, every fixed point for the differential inclusion $x'\in\hbox{Convex Hull}\{f_1(x),f_2(x)\}$ can be approximated by an arbitrarily small two-cycle for the inclusion $x'\in \{f_1(x),f_2(x)\}$, where $f_1,f_2$ are a $C^1$ flows on $\R^n$. 

\end{abstract}
\vfil\break

%%%%%%%%%%%%%%%%%%%%%%%%%%%
%%%%%%%%%%%%%%%%%%%%%%%%%%%
\section{Two-Flows, Stasis Points and Two-Cycles}
%%%%%%%%%%%%%%%%%%%%%%%%%%%
%%%%%%%%%%%%%%%%%%%%%%%%%%%

Consider a controlled system that switches freely between two autonomous flows $x'=f_1(x)$ and $x'=f_2(x)$ where $x\in\R^n$ and $f_i:\R^n\to \R^n$ is $C^1$. This is the differential inclusion 
$$x'\in\{f_1(x),f_2(x)\}.$$
Trajectories are piecewise $C^1$ functions $x(t)$ satisfying this inclusion.  

A \emp{stasis point}is a point $x$ at where the flows are anti-parallel: $$f_1(x)\cdot f_2(x) = -|f_1(x)|\cdot |f_2(x)|.$$ The set of stasis points is the set of possible fixed points for the relaxed differential inclusion $$x'\in\hbox{Convex Hull}(\{f_1(x),f_2(x)\}).$$

A \emp{two-cycle} $x(t)$ is a periodic trajectory $x(t)=x(\delta+t)$ with $0<\kappa<\delta$ and  
$$\begin{array}{rlcc}
	  &x'=f_1(x)&\hbox{for}& 0\le t < \kappa \\
\hbox{and}&x'=f_2(x)&\hbox{for}& \kappa \le t < \delta 
\end{array}
$$

It will be shown that neighborhoods of stasis points generically contain two-cycles.

%%%%%%%%%%%%%%%%%%%%%%%%%%%
%%%%%%%%%%%%%%%%%%%%%%%%%%%
\section{Loops in the Action}
%%%%%%%%%%%%%%%%%%%%%%%%%%%
%%%%%%%%%%%%%%%%%%%%%%%%%%%

An $C^2$ \emp{action} of $\R$ on $\R^n$ is a map $\Phi:\R^n\times\R\to\R^n$ which is $C^2$ in all variables. A \emp{loop} is a point and time $x_0, t_0$ with $\Phi(x_0,t_0)=\Phi(x_0,t_0+\delta)$ for some $\delta>0$. 

A point $x_0$ is \emp{stationary} at time $t_0$ if ${\del\over \del t} \Phi(x_0,t_0)=\vec 0$. Note that a point that remains stationary for some open time interval, ${\del\over \del t} \Phi(x_0,t)\equiv\vec 0$ for $\alpha<t<\beta$, is a (degenerate) loop of all sufficiently small periods. 

The following theorem states that actions with stationary points generically have small loops nearby.  

\begin{theorem}
Let $\Phi:\R^n\times \R \to \R^n$ be a $C^2$ action. Suppose $\exists x_0, t_0$ with ${\partial\over \partial t} \Phi(x_0, t_0)=\vec 0$. If $${\partial^2\over \partial x\partial t} \Phi(x_0, t_0)$$ is non-singular then $\forall \epsilon$, $\exists z_0, \tau_0, \delta>0$ with $\Phi(z_0,\tau_0)=\Phi(z_0,\tau_0+\delta)$ and $$\left|\Phi(z_0,t)-\Phi(x_0,t_0)\right|<\epsilon$$ for $\tau_0\le t<\tau_0+\delta$.
\end{theorem}

\begin{proof}

Without loss of generality we take $x_0=\vec 0$ and $t_0=0$. 

Let 
$$G(x,\delta)= 
\left\{ 
\begin{array}{ccl}
{1\over 2\delta}\left(\Phi(x,\delta)-\Phi(x,-\delta)\right) & \hbox{ if }& \delta\not=0 \cr\cr
{\partial \over \partial t} \Phi (x,0) &  \hbox{ if }& \delta=0
\end{array} 
\right.
$$

So $G$ is $C^1$ at $x=\vec 0,\delta=0$, and ${\partial G \over \partial x}(\vec 0,0)={\partial^2 \Phi \over \partial x\partial t}(\vec 0,0)$ is non-singular by assumption. Therefore the Implicit Function Theorem applies, and since $G(\vec 0,0)=\vec 0$ there exist solutions $x$ to $G(x,\delta)=\vec 0$ for $\delta$ sufficiently small. 

For $\delta\not=0$, a point $x$ with $G(x,\delta)=\vec 0$ has $\Phi(x_0,\delta)=\Phi(x_0,-\delta)$, which is a loop.

\end{proof}

%%%%%%%%%%%%%%%%%%%%%%%%%%%
%%%%%%%%%%%%%%%%%%%%%%%%%%%
\section{Two-Cycle Approximation to Stasis}
%%%%%%%%%%%%%%%%%%%%%%%%%%%
%%%%%%%%%%%%%%%%%%%%%%%%%%%

It is known that relaxed solutions to differential inclusions can be arbitrarily approximated by non-relaxed solutions \cite{fillip,sontag}. Here it is shown that in the case of autonomous two-flows and relaxed fixed points, the approximating trajectory can be taken as a two-cycle. For flows $f_1,f_2$, the idea is to take $\R^n$ modulo the trajectories of $f_1$, and consider the remaining (non-autonomous) action of $f_2$ on $\R^{n-1}$. Loops in this resulting action will pull back to two-cycles in the original system.

\begin{theorem}
Consider two $C^1$ flows $f_1,f_2:\R^n\to\R^n$, and a stasis point $x_0$ with $k_1 f_1(x_0) + k_2 f_2(x_0)=\vec 0$, $k_i > 0$. If ${d \over dx}\left(k_1 f_1 + k_2 f_2\right)$ is non-singular at $x_0$, then there exist two-cycles near $x_0$.
\end{theorem}  

\begin{proof} Without loss of generality we take $x_0= \vec 0$. Since $f_1$ is $C^1$ and non-zero at near $\vec 0$, there is a local diffeomorphism $\rho:\R^n\to\R^n$ with ${d\rho\over dx}\cdot f_1\ident(1,0,\ldots,0)$ near $\vec 0$. Then $g_2={d\rho\over dx}\cdot f_2\circ\rho^{-1}$ is a $C^1$ flow with $g_2(\vec 0)=(-{k_1\over k_2},0,\ldots,0)$. 

We construct the $\R^{n-1}$ action $\Phi$ as follows: given $y= (a_1,\ldots,a_{n-1})\in\R^{n-1}$ let $x(t)=(x_1(t),\ldots,x_n(t))$ be the trajectory under the flow $g_2$ with $x(0) =(0,a_1,\ldots,a_{n-1})\in\R^{n}$. Taking $\Phi(y,t)=(x_2(t),\ldots,x_n(t))$ defines a $C^2$ action on $\R^{n-1}$ with ${\partial\over \partial t} \Phi(\vec 0,0)=\vec 0$. It is straightforward to check that non-singularity of ${d \over dx}\left(\alpha f_1 + \beta f_2\right)$ at $\vec 0\in\R^n$ implies that of ${\partial^2\over \partial x\partial t} \Phi$ at $(\vec 0,0)\in\R^{n-1}\times\R$. 

By Theorem 1, for all $\epsilon>0$ there exist $z=(z_1,\ldots,z_{n-1})$ near $\vec 0\in\R^{n-1}$ with $\Phi(z,\alpha)=\Phi(z,\beta)$ for some $|\alpha|,|\beta|<\epsilon$, $\alpha\not=\beta$. This corresponds to a trajectory $x(t)$ in under $g_2$ on $\R^n$ with $x(0)=(0,z_1,\ldots,z_{n-1})$. Now $x(\alpha)$ and $x(\beta)$ are equal except for their first coordinate, which cannot be zero for sufficiently small $\epsilon$ because $g_2$ is non-zero near $\vec 0$. Thus  $x(\alpha)-x(\beta)$ is parallel to the constant flow $g_1$ and connecting these points with a segment of the $g_1$ flow completes the two-cycle. This maps to a two-cycle for $f_1$ and $f_2$ under $\rho^{-1}$.
\end{proof}

%%%%%%%%%%%%%%%%%%%%%%%%%%%
%%%%%%%%%%%%%%%%%%%%%%%%%%%
\section{Three Questions}
%%%%%%%%%%%%%%%%%%%%%%%%%%%
%%%%%%%%%%%%%%%%%%%%%%%%%%%

For $\R^2$, every two cycle contains a stasis point \cite{johnson}. Is a similar statement true for $\R^n$?

Given a pair of flows $f_1,f_2$, the stasis points generically form one dimensional curves. For $n=2$, an approximating two cycle will cut across this curve twice, and so the set of approximating two-cycles is parameterized by pairs of nearby points on the stasis curve \cite{johnson}. Is there a similar organization of two-cycles in $\R^n$?  

Finally, a natural question is to extend this to $N$-flows. That is, for a differential inclusion $x'\in\{f_1(x),f_2(x),f_3(x)\}$, can a relaxed fixed point always be approximated by a three-cycle?  

%%%%%%%%%%%%%%%%%%%%%%%%%%%%%%%%%%%%%%%%%%%%%%%%%
%%%%%%%%%%%%%%%%%%%%%%%%%%%%%%%%%%%%%%%%%%%%%%%%%
% References
%%%%%%%%%%%%%%%%%%%%%%%%%%%%%%%%%%%%%%%%%%%%%%%%%
%%%%%%%%%%%%%%%%%%%%%%%%%%%%%%%%%%%%%%%%%%%%%%%%%

\end{document}